\documentclass[10pt]{amsart}

\usepackage{amsfonts}

\usepackage{graphicx}
\usepackage{amsmath}
\usepackage{setspace}
\usepackage{enumitem}

\newcommand{\halmos}{{\mbox{\, \vspace{3mm}}} \hfill
\mbox{$\Box$}}

\newtheorem{theorem}{Theorem}
\newtheorem{assumption}{Assumption}

\newtheorem{condition}[theorem]{Condition}
\newtheorem{conjecture}[theorem]{Conjecture}
\newtheorem{corollary}[theorem]{Corollary}

\newtheorem{definition}[theorem]{Definition}
\newtheorem{example}[theorem]{Example}

\newtheorem{lemma}[theorem]{Lemma}

\newtheorem{proposition}[theorem]{Proposition}
\newtheorem{remark}[theorem]{Remark}

\usepackage[english]{babel}
\usepackage[left=3cm,right=3cm,top=2.5cm,bottom=3cm]{geometry}

\usepackage{amsmath,amssymb,amsthm,bbm,color,graphics,version}
\usepackage{mathrsfs}

\newcommand{\bearno}{\begin{eqnarray*}}
\newcommand{\enarno}{\end{eqnarray*}}

\newcommand{\vb}{\vspace{3mm}}

\title{A note on chaotic and predictable representations for It\^{o}-Markov additive processes}

\author{Zbigniew Palmowski}
\address{Faculty of Pure and Applied Mathematics,
Wroc\l{}aw University of Science and Technology,
Wyb. Wyspia\'nskiego 27, 50-370 Wroc\l{}aw, Poland}
\email{zbigniew.palmowski@gmail.com}

\author{{\L}ukasz Stettner}
\address{Institute of Mathematics Polish Acad. Sci., Sniadeckich 8, 00-656 Warsaw, also Vistula University}
\email{stettner@impan.pl}

\author{Anna Sulima}
\address{Institute of Mathematics, Jagiellonian University, \L{}ojasiewicza 6, 30-348 Cracow and Institute of Economics, Polish Academy of Sciences, Parade Square 1, 00-901 Warsaw}
\email{anna.sulima@wp.pl}

\thanks{
This work was partially supported by the National Science Centre under the grant\linebreak
 2015/17/B/ST1/01102.
}

\date{}

\begin{document}

\begin{abstract}
In this paper we provide predictable and chaotic representations
for It\^{o}-Markov additive processes $X$.
Such a process is governed by a finite-state CTMC $J$ which allows one to modify the parameters of the It\^{o}-jump process (in so-called regime switching manner). In addition, the transition of $J$ triggers the jump
of $X$ distributed depending on the states of $J$ just prior to the transition.
This family of processes includes Markov modulated It\^{o}-L\'evy processes and Markov additive processes.
The derived chaotic representation of a square-integrable random variable is given as a sum of
stochastic integrals with respect to some explicitly constructed orthogonal martingales.
We identify the predictable representation of a square-integrable martingale
as a sum of stochastic integrals of predictable processes with respect to
Brownian motion and power-jumps martingales related to all the jumps
appearing in the model. This result generalizes the seminal result of Jacod-Yor and is of importance in financial mathematics. The derived representation then allows one to enlarge the incomplete market by a series of power-jump
assets and to price all market-derivatives.
\vspace{3mm}

\noindent


\end{abstract}

\maketitle

\pagestyle{myheadings} \markboth{\sc Z.\ Palmowski --- {\L}.\ Stettner --- A.\ Sulima} {\sc
Chaotic and predictable representations for It\^{o}-Markov additive processes}

\vspace{1.8cm}

\tableofcontents

\newpage

\section{Introduction}

\text{ } \text{ } The Martingale Representation Theorem allows one to represent all martingales, in a certain space, as stochastic integrals.
This result is the cornerstone of mathematical finance that produces hedging of some derivatives (see e.g. F\"ollmer and Schied \cite{FS}).
It is also the basis for the theory of backward stochastic differential equations (see El Karoui et al. \cite{EPQ} for a review).
The martingale representation theorem can also lead to Wiener-like chaos expansions, which appear in the Malliavin calculus (
see Di Nunno et al. \cite{DOP}).

\text{ } \text{ } The martingale representation theorem is commonly considered in a Brownian setting.
It states that any square-integrable martingale $M$ adapted to the natural filtration of the Brownian motion $W$
can be expressed as a stochastic integral with respect to the same Brownian motion $W$:
\begin{equation}\label{86}
M(t)=M(0)+\int_{0}^{t}\xi(s) \mathrm{d}W(s),
\end{equation}
where $\xi$ satisfies suitable integrability conditions.
This is the classical It\^{o} theorem (see Rogers and Williams \cite[Thm. 36.1]{Rogers}).

\text{ } \text{ } For processes having jumps coming according to a compensated Poisson measure $\bar{\Pi}$,
the representation theorem was originally established by Kunita
and Watanabe in their seminal paper \cite[Prop. 5.2]{KW} (see also Applebaum \cite[Thm. 5.3.5]{App}).

\text{ } \text{ } In particular, on probability space $(\Omega, \mathcal{F}, \mathbb{P})$ with augmented natural filtration of a L\'evy process $\{\mathcal{H}_{t}\}_{t\geq 0}$, any locally square-integrable $\{\mathcal{H}_{t}\}$-martingale $M$ has the following representation:
\begin{equation}\label{87}
M(t)=M(0)+\int_{0}^{t}\xi_{1}(s) \mathrm{d} W(s)+\int_{0}^{t}\int_{\mathbb{R} \setminus \{0\}}\xi_{2}(s,x) \bar{\Pi}(\mathrm{d}s,\mathrm{d}x),
\end{equation}
where  $\xi_{1}$ and $\xi_{2}$ are $\{\mathcal{H}_{t}\}$- predictable processes satisfying suitable integrability conditions.
Later the results were generalized in a series of papers by Chou and Meyer \cite{CM}, Davis \cite{Da} and Elliott \cite{E1, E2}. In fact, Elliott \cite{E1} was assuming that there are a finite number of jumps within each finite interval.
To deal with a general L\'evy process (of unbounded variation) Nualart and Schoutens \cite{NS} derived a martingale representation
using the chaotic representation property
in terms of a suitable orthogonal sequence of martingales, obtained as the orthogonalization of the compensated
Teugels power jump processes of the L\'evy process (see also Corcuera et al. \cite{CNS}).
More on martingale representation can be found in Jacod and Shiryaev \cite{JS} , Liptser and Shiryaev \cite{LS}, Protter \cite{P}
and in a review by Davis \cite{D}.

\vb
\text{ } \text{ } The goal of this paper is to derive the martingale representation theorem for an It\^{o}-Markov additive process $X$.  It\^{o}-Markov additive processes are a natural generalization of
It\^{o}-L\'evy processes and hence also of L\'evy processes. The use of It\^{o}-Markov additive processes is widespread, making them a classical model in applied probability with a variety
of application areas, such as queues, insurance risk, inventories, data communication, finance,
environmental problems and so forth (see Asmussen \cite{A2}, Asmussen and  Albrecher \cite{A1}, Asmussen et al. \cite{asm_avram_pist},
Asmussen and Kella \cite{AK}, \c Cinlar \cite{Ca, Cb}, Prabhu \cite[Ch. 7]{prabhu},
Pacheco and Prabhu \cite{PP}, Palmowski and Rolski \cite{PR}, Palmowski and Ivanovs \cite{meiv} and references therein).

\text{ } \text{ }In this paper, we will also derive a chaotic representation of any square-integrable random variable in terms of certain
orthogonal martingales. It should be underlined that the derived representation is of predictable type,
that is, any square-integrable martingale is a sum of stochastic integrals of some predictable processes with respect to explicit martingales. This representation is hence different from the one given in
\eqref{87} (in the jump part) and it can be used to prove completeness of some financial markets.

\text{ } \text{ } This paper is organized as follows.
In Section \ref{prel} we present the main results.
Their proofs are given in Section \ref{sec:proos}, preceeded by a crucial preliminary fact
presented in Section \ref{pol3}.

\section{Main results}\label{prel}
\text{ } \text{ } We begin by defining It\^{o}-Markov additive processes. Let $(\Omega, \mathcal{F}, \mathbb{P})$ be a complete probability space and let $\mathbb{T} := [0;T]$ be a finite time horizon, where $0<T<\infty$ is fixed. On this probability space we consider a homogeneous continuous-time Markov chain  $J :=\{J(t): t\in \mathbb{T}\} $ with a finite state space. For simplicity, we follow the notation of Elliott et al. \cite{EAM} and we identify the state space as a finite set of unit vectors $E := \{\textbf {e}_{1}, \ldots, \textbf {e}_{N}\}$ of $\mathbb{R}^{N}$. Here the $j$th component of $\textbf{e}_{i}$ is the Kronecker delta $\delta_{ij}$ for each $i,j = 1, 2, \ldots ,N$. Moreover, the Markov chain  $J$ is characterized by an intensity matrix $[\lambda_{ij}]_{i,j=1,2, \ldots ,N}$. The element $\lambda_{ij} $ is the transition intensity of the Markov chain $J$ jumping from state $\textbf {e}_{i}$ to state $\textbf {e}_{j}$. \\
A process $(J,X)=\{(J(t),X(t)): t\in \mathbb{T}\}$ on the state space $\{\textbf{e}_1,\ldots, \textbf{e}_N\} \times \mathbb{R}$ is a {\it Markov additive process} (MAP) if $(J,X)$ is a Markov process and the conditional distribution of $(J(s+t),X(s+t)-X(s))$ for $s, t\in \mathbb{T}$, given $(J(s),X(s))$, depends only on $J(s)$ (see \c Cinlar \cite{Ca, Cb}).   
Every MAP has a very special structure. It is usually said that $X$ is the additive component and $J$ is the background process representing the environment.
Moreover, the process $X$ evolves as a L\'{e}vy process while $J(t) = \textbf{e}_j$.
Following Asmussen and Kella \cite{AK} we can decompose the process $X$ as follows:
\begin{equation}\label{defX}
X(t)=\overline{X}(t)+\overline{\overline{X}}(t),
\end{equation}
where
\begin{equation}\label{ovvX}
\overline{\overline{X}}(t):=\sum_{i=1}^{N}\Psi_{i}(t)
\end{equation}
for
\begin{equation}\label{Psij}
\Psi_{i}(t):=\sum_{n\geq1}U^{(i)}_{n}\textbf{1}_{\{J(T_{n})=\textbf {e}_{i},\textrm{ } T_{n}\leq t\}}
\end{equation}
and for the jump epochs $\{T_n\}$ of $J$.
Here $U^{(i)}_{n}$ $(n\geq 1, 1\leq i \leq N)$ are independent random variables which are also independent of $\overline{X}$ such that for every $i$, the random variables $U_{n}^{(i)}$ are identically distributed. 
Note that we can express the process $\Psi_{i}$ as follows:
\begin{equation*}
\Psi_{i}(t)=\int_0^t \int_{\mathbb{R}}x \textrm{ } \Pi_{U}^{i} (\mathrm{d}s, \mathrm{d}x)
\end{equation*}
for the point measures
\begin{equation}\label{PiU}
\Pi_{U}^{i} ([0,t], \mathrm{d}x):=\sum_{n\geq1}  \mathbb{P}(U^{(i)}_{n} \in \mathrm{d}x )  \textbf{1}_{\{J(T_{n})=\textbf {e}_{i},\textrm{ } T_{n}\leq t\}}, \textrm{ }\textrm{ } i = 1, \ldots ,N.
\end{equation}

\begin{remark}
One can consider
jumps $U^{(ij)}$ with distribution depending also on the state $\textbf{e}_j$ the Markov chain is jumping to by
extending the state space to the pairs $(\textbf{e}_i, \textbf{e}_j)$ (see Gautam et al. \cite[Thm. 5]{Palmowskietal}
for details).
\end{remark}

\text{ } \text{ } The first component in equation (\ref{defX}) is an It\^{o}-L\'evy process and it has the following decomposition (see Oksendal and Sulem \cite[p. 5]{OS}):
\begin{equation} \label{1010}
\overline{X}(t):= \overline{X}(0)+ \int_{0}^{t} \mu_{0}(s)\mathrm{d}s + \int_{0}^{t}\sigma_{0}(s)\mathrm{d}W(s)+ \int_{0}^{t} \int_{\mathbb{R}} \gamma(s-,x) \bar{\Pi}(\mathrm{d}s,\mathrm{d}x),
\end{equation}
where $W$ denotes the standard Brownian motion independent of $J$, $ \bar{\Pi}(\mathrm{d}t,\mathrm{d}x):=\Pi(\mathrm{d}t,\mathrm{d}x)-\nu(\mathrm{d}x)\mathrm{d}t$ is the compensated Poisson random measure which is independent of $J$ and $W$ and
\begin{equation}\label{m0}
 \mu_{0}(t):=\langle  \boldsymbol\mu_{0}, J(t) \rangle=\sum_{i=1}^{N}\mu_{0}^{i}\langle \textbf{e}_{i},J(t)\rangle,
 \end{equation}
\begin{equation}\label{sigma0}
\sigma_{0}(t):=\langle \boldsymbol\sigma_{0}, J(t) \rangle =\sum_{i=1}^{N}\sigma_{0}^{i} \langle \textbf{e}_{i},J(t) \rangle,
\end{equation}
\begin{equation}\label{gamma0}
\gamma(t,x):=\langle  \boldsymbol\gamma(x), J(t) \rangle=\sum_{i=1}^{N}\gamma_{i}(x) \langle \textbf{e}_{i}, J(t) \rangle,
\end{equation}
for some vectors
$\boldsymbol\mu_{0}=(\mu_{0}^{1}, \ldots ,\mu_{0}^{N})'\in \mathbb{R}^{N}$, $\boldsymbol\sigma_{0}=(\sigma_{0}^{1}, \ldots , \sigma_{0}^{N})'\in\mathbb{R}^{N}$ and the vector-valued measurable function $\boldsymbol\gamma(x)=(\gamma_{1}(x), \ldots ,\gamma_{N}(x))$.
The measure $\nu$ is the so-called jump-measure identifying the distribution of the sizes of the jumps of the Poisson measure $\Pi$.
The components $\overline{X}$ and $\overline{\overline{X}}$ in \eqref{defX} are, conditionally on the state of the Markov chain  $J$, independent.
Additionally, we suppose that the L\'{e}vy measure satisfies, for some $ \varepsilon >0$ and $\lambda >0$,
\begin{equation}\label{42}
\int_{(- \varepsilon,  \varepsilon)^{c}} \exp(\lambda |\gamma(s-,x)|)\nu(\mathrm{d}x)<\infty, \mathbb{P}-a.e., \quad \quad \quad \int_{(- \varepsilon,  \varepsilon)^{c}} \exp(\lambda x)\mathbb{P}(U^{(i)}\in \mathrm{d}x)<\infty 
\end{equation}
for $i=1,\ldots,N$, where $U^{(i)}$ is a generic size of jump when the Markov chain $J$ is jumping into state $i$.
This implies that
\begin{equation*}
 \int_{\mathbb{R}}|\gamma(s-,x)|^{k}\nu( \mathrm{d}x)<\infty, \mathbb{P}-a.e.,\quad \quad \quad \mathbb{E(}U^{(i)})^k<\infty, \textrm{ }\textrm{ } \textrm{ }\textrm{ }i=1,\ldots,N, \textrm{ }\textrm{ }  k \geq 2
\end{equation*}
and that the characteristic function $\mathbb{E}[\exp(ku X]$ is analytic in a neighborhood of $0$. Moreover, $X$
has moments of all orders and the polynomials are dense in $L^{2}(\mathbb{R}, \mathrm{d} \varphi (t,x))$, where $\varphi(t,x):=\mathbb{P}(X(t)\leq x)$.\\

\text{ } \text{ } Now we are ready to define the main process in this paper.
The process  $(J,X)= \{(J(t),X(t)):t \in \mathbb{T} \}$ (for simplicity sometimes we write only $X$) with the decomposition (\ref{defX}) is called an {\it It\^{o}-Markov additive process}.

\text{ } \text{ } This process evolves as the It\^{o}-L\'evy process $\overline{X}$ between changes of states of the Markov chain $J$, that is, its parameters depend on the current state $\textbf{e}_i$ of the Markov chain $J$.
In addition, a transition of $J$ from $\textbf{e}_i$ to $\textbf{e}_j$ triggers a jump of $X$ distributed as $U^{(i)}$. This is a so-called non-anticipative It\^{o}-Markov additive process. 

\text{ } \text{ } It\^{o}-Markov additive processes are a natural generalization of It\^{o}-L\'evy processes and thus of L\'evy processes.
Moreover, if $\gamma(s,x)=x$ then $X$ is a Markov additive process. If additionally $N=1$, then $X$ is a L\'evy process.
If $U^{(j)}\equiv 0$ and $N> 1$ then $X$ is a Markov modulated L\'{e}vy process (see Pacheco et al. \cite{PTP}).
If there are no jumps additionally, that is, $ \bar{\Pi}( \mathrm{d}s, \mathrm{d}x)=0$, we have a Markov modulated Brownian motion.

\text{ } \text{ } From now, we will work with the following filtration on $(\Omega, \mathcal{F}, \mathbb{P})$:
\begin{equation}
\mathcal{F}_{t}:= \mathcal{G}_{t} \vee  \mathcal{N},
\end{equation}
where $ \mathcal{N}$ are the $\mathbb{P}$-null sets of $\mathcal{F}$ and
\begin{equation}
\mathcal{G}_{t}:=\sigma\{J(s),W(s), \Gamma(s), \Pi_{U}^{1} ([0,s], \mathrm{d}x), \ldots, \Pi_{U}^{N} ([0,s], \mathrm{d}x) ;s\leq t\}
\end{equation}
for
\begin{equation}
\Gamma(t):= \int_{0}^{t} \int_{\mathbb{R}} \gamma(s-,x) \Pi(\mathrm{d}s,\mathrm{d}x).
\end{equation}
Note that the filtration $ \{\mathcal{F}_{t}\}_{ t\geq 0 }$ is right-continuous (see Karatzas and Shreve \cite[Prop. 7.7 ]{KS} and also Protter \cite[Thm. 31 ]{P}). By the same arguments as in the proof of Thm. 3.3 of Liao \cite{Liao},
the filtration $ \{\mathcal{G}_{t}\}_{ t\geq 0 }$  is quivalent to
\begin{equation}
\mathcal{G}_{t}:=\sigma\{J(s), \overline{X}(s), \Pi_{U}^{1} ([0,s], \mathrm{d}x), \ldots, \Pi_{U}^{N} ([0,s], \mathrm{d}x) ;s\leq t\}.
\end{equation}

\text{ } \text{ } To present the main result we need a few additional processes.
We observe that the Markov chain $J$ can be represented in terms of a marked point process $\Phi_{j} $ defined by
\begin{equation*}
\Phi_{j}(t):=\Phi([0,t] \times \textbf e_{j})=\sum_{n \geq 1}\textbf{1}_{\{J(T_{n}) =  \textbf {e}_{j}, \textrm{ }   T_{n} \leq t \}}=\Pi_{U}^{j} ([0,t], \mathbb{R}), \textrm{ }\textrm{ } j = 1,2, \ldots ,N
\end{equation*}
for the jump epochs $\{T_n\}$ of $J$.
Note that the process $\Phi_{j} $ describes the number of jumps into state $\textbf{e}_{j}$ up to time $t$.
Let $\phi_{j}$ be the dual predictable projection of $\Phi_{j}$ (sometimes called the compensator). That is, the process
\begin{equation}\label{ovlPhi}
\overline{\Phi}_{j}(t):=\Phi_{j}(t)-\phi_{j}(t),\textrm{ }\textrm{ } j = 1, \ldots ,N,
\end{equation}
is an $\{\mathcal{F}_{t}\}$-martingale and it is called the $j$th Markovian jump martingale. Note that $\phi_{j} $ is unique and
\begin{equation}\label{9878}
\phi_{j}(t):= \int_{0}^{t}\lambda_{j}(s)\mathrm{d}s
\end{equation}
for 
\begin{equation}\label{lamb}
\lambda_{j}(t):=\sum\limits_{i\neq j} \textbf{1}_{\{J(t-)=  \textbf{e}_{i}\}}\lambda_{ij} 
\end{equation}
(see Zhang et al. \cite[p. 290]{ZESG}).\\
Following Corcuera et al. \cite{CNS} we introduce the power-jump processes
$$X^{(k)}(t):= \sum\limits_{0<s\leq t}(\Delta \overline{X}(s))^{k},  \textrm{ } \textrm{ } \textrm{ }\textrm{ }\textrm{ }\textrm{ } k \geq 2,$$
where $\Delta \overline{X}(s)=\overline{X}(s)-\overline{X}(s-)$. We set $X^{(1)}(s)=\overline{X}(s)$.
The process $X^{(k)}$ is also an It\^{o}-L\'evy process with the same jump times as the original process $\overline{X}$ but with their sizes being
the $k$th powers of the jump sizes of $\overline{X}$. From Protter \cite[p. 29]{P} we have

\begin{equation*}
\mathbb{E}\big[X^{(k)}(t)\big|  \mathcal{J}_{t} \big]=\mathbb{E}\bigg( \sum_{0<s\leq t} (\Delta \overline{X}(s))^{k} \big|  \mathcal{J}_{t}  \bigg)=
\int_0^t \int_{\mathbb{R}}\gamma^{k}(s-,x)\nu(\mathrm{d}x)\mathrm{d}s <\infty, \mathbb{P}-a.e. \textrm{ } \textrm{ } \textrm{ }\textrm{ }\textrm{ }\textrm{ } k \geq 2,
\end{equation*}

for $\mathcal{J}_{t}:=\sigma\{ J(s): s \leq t\}$ and hence the processes
\begin{equation}\label{ovlX}
\overline{X}^{(k)}(t):=X^{(k)}(t)-\int_0^t \int_{\mathbb{R}}\gamma^{k}(s-,x)\nu(\mathrm{d}x)\mathrm{d}s, \textrm{ } \textrm{ } \textrm{ }\textrm{ }\textrm{ }\textrm{ } k \geq 2,
\end{equation}
are $\{\mathcal{F}_{t}\}$-martingales (called Teugels martingales of order $k$; see Schoutens \cite{S} for details).
Indeed, since $\overline{X}$ is $\{\mathcal{F}_{t}\}$-adapted by Jacod and Shiryaev \cite[Prop. 1.25]{JS} the process  $\Delta \overline{X}$ is also
$\{\mathcal{F}_{t}\}$-adapted. Furthermore, any integral of an $\{\mathcal{F}_{t}\}$-adapted process with respect to an $\{\mathcal{F}_{t}\}$-adapted stochastic measure or an $\{\mathcal{F}_{t}\}$-adapted stochastic process is still adapted (see Jacod and Shiryaev  \cite[Prop. 3.5 and Thm 4.31(i)]{JS}). Hence $X^{(k)}$ is $\{\mathcal{F}_{t}\}$-adapted for $k \geq 2$.\\
We will also need power martingales related to the second component of $X$ given in \eqref{defX}, namely to $\overline{\overline{X}}$
or to $\Psi_i$ defined in \eqref{Psij}. For $l \geq 1$ and $i=1, \ldots , N$ we define
\begin{equation*}\label{Psijk}
\Psi_{i}^{(l)}(t): =\sum_{n\geq1}\left(U^{(i)}_{n}\right)^l\textbf{1}_{\{J(T_{n})=\textbf {e}_{i},\textrm{ } T_{n}\leq t\}}=\int_0^t \int_{\mathbb{R}}x^{l} \textrm{ } \Pi_{U}^{i} (\mathrm{d}s, \mathrm{d}x)
\end{equation*}
for $\Pi_{U}^{i}$ given by (\ref{PiU}).
The compensated version of $\Psi_{i}^{(l)}$ is called an impulse regime-switching martingale if 
\begin{equation}\label{ovlPsi}
\overline{\Psi}_{i}^{(l)}(t):= \Psi_{i}^{(l)}(t) - \mathbb{E}\big(U^{(i)}_{n}\big)^{l}\phi_{i}(t)=\int_0^t \int_{\mathbb{R}}x^{l} \textrm{ }\bar{ \Pi}_{U}^{i} (\mathrm{d}s, \mathrm{d}x),
\end{equation}
where $\bar{ \Pi}_{U}^{i} (\mathrm{d}t, \mathrm{d}x) = \Pi_{U}^{i}(\mathrm{d}t, \mathrm{d}x) - \lambda_{i}(t)\eta_{i}(\mathrm{d}x)\mathrm{d}t  $  for $\lambda_{i}$ defined in (\ref{lamb}) and $\eta_{i}(\mathrm{d}x)= \mathbb{P}(U^{(i)}_{n} \in \mathrm{d}x ) $. 

\text{ } \text{ } Using similar arguments to those above it follows that $\overline{\Psi}_{i}^{(l)}$ is an $\{\mathcal{F}_{t}\}$-martingale for $l \geq 1$ and $i=1,...,N$.\\
\text{ } \text{ }Corcuera et al. \cite{CNS} motivate trading in power-jump assets as follows. 
Power-jump process of order two is just a variation process of degree two, i.e. a quadratic variation process (see Barndorff-Nielsen and Shephard \cite{BNSh, BNSh2} ), and is related to the so-called realized variance. Contracts on realized variance have found their way into OTC markets and are now traded regularly. Typically a 3th-power-jump asset measures a kind of asymmetry ("skewness") and a 4th-power-jump process measures extremal movements ("kurtosis"). Trade in such assets can be of use if one likes to bet on the realized skewness or realized kurtosis of the stock. Furthermore, an insurance contract against a crash can also be  easily built from 4th-power jump (or ith-power jump, $i > 4$) assets.
\vb

\text{ } \text{ } We denote by $\mathcal{M}^2$ the set of square-integrable $\{\mathcal{F}_{t}\}$-martingales, i.e. $M \in \mathcal{M}^2$  if $M$ is a martingale, $M(0)=0$ and $\sup_t \mathbb{E} M^2 (t)<\infty$. 
The martingale convergence theorem implies that each $M \in \mathcal{M}^2$ is closed, i.e. there is an $\{\mathcal{F}_{t}\}$-measurable random variable $M$  such that $M(t) \rightarrow M(\infty)$ in $L^2 (\Omega, \mathcal{F})$ and for each $t$, $M(t)=\mathbb{E}[M(\infty)|\mathcal{F}_t]$.
Thus there is a one-to-one correspondence between $\mathcal{M}^2$  and $L^2 (\Omega, \mathcal{F})$, so that $\mathcal{M}^2$ is a Hilbert space under the inner product $M_1 (t) \cdot M_2 (t)=\mathbb{E}[M_1(\infty)M_2(\infty)]$.
Following Protter \cite[p. 179]{P}, we say that two martingales $M_1, M_2 \in \mathcal{M}^2$ are \emph{strongly orthogonal} if their product $M_1\cdot M_2$ is a uniformly integrable martingale. As noted in Protter \cite{P}, $M_1, M_2\in \mathcal{M}^2$ are strongly orthogonal if and only if $[M_1, M_2]$ is a uniformly integrable martingale. We say that two random variables $Y_1, Y_2 \in L^2 (\Omega, \mathcal{F})$ are \emph{weakly orthogonal} if $\mathbb{E} [ Y_1(t), Y_2(t)]=0$. Clearly, strong orthogonality implies weak orthogonality.
\\
One can obtain a set $ \{H^{(k)}, k \geq 1\}$ of pairwise strongly orthonormal martingales such that each
$ H^{(k)}$ is a linear combination of the  $\overline{X}^{(n)} (n=1,\ldots,k)$,
\begin{equation} \label{defH}
H^{(k)}=a_{k,k}\overline{X}^{(k)}+a_{k,k-1}\overline{X}^{(k-1)}+\cdots+a_{k,1}\overline{X}^{(1)},\textrm{ } \textrm{ } \textrm{ }\textrm{ }\textrm{ }\textrm{ } k \geq 1.
\end{equation}
The constants $a_{\cdot, \cdot}$ can be calculated as described in  Schoutens \cite{S} -- they correspond to the coefficients
of the orthonormalization of the polynomials $1,x,x^{2}, \ldots $.
Similarly we  proceed with the martingales $\overline{\Psi}_{i}^{(l)}$ and $\overline{\Phi}_{i}$ and
construct pairwise strongly orthonormal martingales $G_{i}^{(l)}$ ($i=1,\ldots,N$, $l \geq 1$) that are appropriate
linear combination of the processes $\overline{\Psi}_{i}^{(l)}$ and $\overline{\Phi}_{i}$:
\begin{equation} \label{defG}
G_{i}^{(l)}=c^{(i)}_{l,l}\overline{\Psi}_{i}^{(l-1)}+c^{(i)}_{l,l-1}\overline{\Psi}_{i}^{(l-2)}+\cdots+c^{(i)}_{l,2}\overline{\Psi}_{i}^{(1)}+c^{(i)}_{l,1} \overline{\Phi}_{i}.
\end{equation}
We will find the coefficients $c_{\cdot,\cdot}^{(i)}$ as follows.
Fix $i\in\{1,\ldots, N\}$.
Let us consider two spaces. The first one is the space $S_{1}$ of all real polynomials
on the positive real line endowed with the scalar product 
\begin{equation*}
\langle P(x), Q(x) \rangle_{1}:= \mathbb{E}\big(P(U^{(i)})Q(U^{(i)})\big)\mathbb{E} \big(\Phi_{i}(1)\big).
\end{equation*}
Note that
\begin{equation*}
\langle x^{l}, x^{h} \rangle_{1}= \mathbb{E}(U^{(i)})^{l+h}\mathbb{E} \big(\Phi_{i}(1)\big).
\end{equation*}
The other space, $S_{2}$, is the space of all linear transformations of the processes  $\overline{\Psi}_{i}^{(l)}$ and $\overline{\Phi}_{j}$, i.e.,
\begin{equation*}
S_{2}= \{c_{l}\overline{\Psi}_{i}^{(l-1)}+c_{l-1}\overline{\Psi}_{i}^{(l-2)}+\cdots+c_{2}\overline{\Psi}_{i}^{(1)}+c_{1} \overline{\Phi}_{i};\quad l\in \{ 1, 2, \ldots  \},\quad c_{i} \in \mathbb{R}\}.
\end{equation*}
We endow this space with the scalar product 
\begin{align*}
\langle \overline{\Psi}_{i}^{(l)}(t), \overline{\Psi}_{i}^{(h)}(t) \rangle_{2} &:= \mathbb{E}\Big([\overline{\Psi}_{i}^{(l)}, \overline{\Psi}_{i}^{(h)}  ](1)\Big)=\mathbb{E}(U^{(i)})^{l+h}\mathbb{E} \big(\Phi_{i}(1)\big),\\
\langle \overline{\Psi}_{i}^{(l)}(t),  \overline{\Phi}_{i}(t) \rangle_{2} &:= \mathbb{E}\Big([\overline{\Psi}_{i}^{(l)}, \overline{\Phi}_{i}  ](1)\Big)= \mathbb{E}(U^{(i)})^{l} \mathbb{E} \big(\Phi_{i}(1)\big),\\
\langle \overline{\Phi}_{i}(t),  \overline{\Phi}_{i}(t) \rangle_{2}&:= \mathbb{E}\Big([\overline{\Phi}_{i}, \overline{\Phi}_{i}  ](1)\Big)=\mathbb{E}\big(\Phi_{i}(1)\big),
\end{align*}
for $i=1,\ldots,N$ and $l,h \geq 0 $.
One clearly sees that $x^{l}\leftrightarrow\overline{\Psi}_{i}^{(l)} $ is an isometry between $S_{1}$ and $S_{2}$. An orthogonalization of $\{1,x,x^{2},\ldots\}$ in $S_{1}$ gives an orthogonalization of $\{\overline{\Phi}_{i}, \overline{\Psi}_{i}^{(1)},\overline{\Psi}_{i}^{(2)},\ldots\}$.\\
\text{ } \text{ } Finally, for $i \neq j$, the processes $\overline{\Psi}_{i}^{(l)} $, $\overline{\Phi}_{i}$ and $\overline{\Psi}_{j}^{(h)}, \overline{\Phi}_{j}$ do not jump at the same time, so
\[[\overline{\Psi}_{i}^{(l)}, \overline{\Psi}_{j}^{(h)}](t)=\sum_{s\leq t} \Delta \overline{\Psi}^{(l)}_{i}(s)\Delta \overline{\Psi}^{(h)}_{j}(s)=0,  \quad \quad \quad [\overline{\Psi}_{i}^{(l)}, \overline{\Phi}_{j}](t)=\sum_{s\leq t} \Delta \overline{\Psi}^{(l)}_{i}(s) \Delta \overline{\Phi}_{j}(s)=0 \] 
and $$[\overline{\Phi}_{i}, \overline{\Phi}_{j}](t)=\sum\limits_{s\leq t} \Delta \overline{\Phi}_{i}(s) \Delta \overline{\Phi}_{j}(s)=0.$$ For the same reason we have
\[[\overline{\Psi}_{i}^{(l)}, H^{(k)}](t)=0, \text{ } [\overline{\Phi}_{i}, H^{(k)}](t)=0\qquad\text{for $i\in\{1,\ldots, N\}$, $k \geq 1$ and $l \geq 1 $}.\]
\text{ } \text{ } In this way all martingales in \eqref{defH} and \eqref{defG} are pairwise strongly orthogonal.\\

\text{ } \text{ } The main results of this paper is given in the next two theorems.
\begin{theorem}\label{987}
Any square-integrable $\left\{ \mathcal{F}_t \right\}$-measurable random variable $F$ can be represented as follows:
\begin{eqnarray}\label{22a}
\nonumber F(t) &=& \mathbb{E}[F(t)]+ \sum_{i=1}^{N} \sum_{s=1}^{\infty} \sum_{\tau=1}^{\infty}  \sum_{\iota_{1},\ldots,\iota_{s}\geq 1} \sum_{\upsilon_{1},\ldots,\upsilon_{\tau}\geq 1} \int_{0}^{t} \int_{0}^{t_{1}-}\ldots \int_{0}^{t_{s+\tau-1-}}f_{(\upsilon_{1},\ldots,\upsilon_{\tau}, \iota_{1},\ldots,\iota_{s}, i )}(t_{1},t_{2},\ldots,t_{s+\tau })\\
& & \qquad\qquad \mathrm{d}G^{(\upsilon_{\tau})}_{i}(t_{s+\tau})\ldots\mathrm{d}G^{(\upsilon_{1})}_{i}(t_{s+1})\mathrm{d}H^{(\iota_{s})}(t_{s}) \ldots\mathrm{d}H^{(\iota_{1})}(t_{1}),
\end{eqnarray}
where  $f_{(\upsilon_{1},\ldots,\upsilon_{\tau}, \iota_{1},\ldots,\iota_{s} )}$ are some random fields for which (\ref{22a}) is well-defined on $L^{2}(\Omega,\mathcal{F})$, processes $G^{(l)}_{i}$ and $H^{(k)}$
($i=1,\dots, N$, $l,k\geq 1$) are orthogonal martingales and the convergence is in $L^2$ sence.
\end{theorem}

\begin{remark}\label{rem34}\rm
The right-hand side of (\ref{22a}) is understood as follows. We take a finite sum
\begin{eqnarray}\label{22abc}
& & \sum_{i=1}^{N} \sum_{s=1}^{A_1} \sum_{\tau=1}^{A_2}  \sum_{\iota_{1},\ldots,\iota_{s}\geq 1
} \sum_{\upsilon_{1},\ldots,\upsilon_{\tau}\geq 1} \int_{0}^{t} \int_{0}^{t_{1}-} \cdots \int_{0}^{t_{s+\tau-1-}}f_{(\upsilon_{1},\ldots,\upsilon_{\tau}, \iota_{1},\ldots,\iota_{s}, i)}(t_{1},t_{2},\ldots,t_{s+\tau })\nonumber \\
\nonumber & & \qquad\qquad \mathrm{d}G^{(\upsilon_{\tau})}_{i}(t_{s+\tau}) \ldots \mathrm{d}G^{(\upsilon_{1})}_{i}(t_{s+1})\mathrm{d}H^{(\iota_{s})}(t_{s}) \ldots \mathrm{d}H^{(\iota_{1})}(t_{1})
\end{eqnarray}
in $L^2 (\Omega, \mathcal{F})$. Since $L^2 (\Omega, \mathcal{F})$ is a Hilbert space, the right-hand side of  (\ref{22a}) is understood as the limit of the above expression in $L^2 (\Omega, \mathcal{F})$ for $A_1 \rightarrow \infty$ and $A_2 \rightarrow \infty $.
\end{remark}

\begin{theorem}\label{8}
Any square-integrable $\left\{ \mathcal{F}_t \right\}$-martingale $M$ can be represented as follows:
\begin{eqnarray}\label{b8}
\nonumber  M(t) &=& M(0)+ \int_{0}^{t} h_{X}^{(1)}(s)\mathrm{d}X(s)+\sum_{k=2}^{\infty}\int_{0}^{t}h^{(k)}_{X}(s)\mathrm{d}\overline{X}^{(k)}(s) +\sum_{j=1}^{N}\int_{0}^{t}h_{\Phi}^{(j)}(s)\mathrm{d}\overline{\Phi}_{j}(s)\\
& & +\sum_{l=1}^{\infty}\sum_{i=1}^{N}\int_{0}^{t} h^{(l,i)}_{\Psi}(s)\mathrm{d}\overline{\Psi}_{i}^{(l)}(s),
\end{eqnarray}
where $h_{X}^{(1)}$, $h_{X}^{(k)}$, $h_{\Phi}^{(j)}$ and $h^{(l,i)}_{\Psi}$ (for $i,j=1,\ldots,N $, $k \geq 2$ and $l \geq 1$) are predictable processes.
\end{theorem}

\begin{remark}\label{rem3} \rm
The right-hand side of (\ref{b8}) is understood in the same way as in Remark \ref{rem34}, that is, it is the limit in $\mathcal{M}^2$ of
\begin{eqnarray*}
 \int_{0}^{t} h_{X}^{(1)}(s)\mathrm{d}\overline{X}(s)+\sum_{k=2}^{B_1}\int_{0}^{t}h^{(k)}_{X}(s)\mathrm{d}\overline{X}^{(k)}(s) +\sum_{j=1}^{N}\int_{0}^{t}h_{\Phi}^{(j)}(s)\mathrm{d}\overline{\Phi}_{j}(s) +\sum_{l=1}^{B_2}\sum_{i=1}^{N}\int_{0}^{t} h^{(l,i)}_{\Psi}(s)\mathrm{d}\overline{\Psi}_{i}^{(l)}(s)
\end{eqnarray*}
for $B_1 \rightarrow \infty$ and $B_2 \rightarrow \infty $.
\end{remark}

\text{ } \text{ } The proof of this theorem will be given in Section \ref{pol1}.
The main idea of the proof is that every random variable $F$ in $L^{2}(\Omega,\mathcal{F}_{t})$ can be approximated by some type of polynomials. 
For these polynomials we will use the It\^{o} formula together with induction to get an appropriate representation first
in terms of orthogonal multiple stochastic integrals and then as a sum of single stochastic integrals.
The main step of the construction is given in Section \ref{pol3}.

\begin{remark}\rm
Note that the representation in Theorem \ref{8} is as a sum of integrals with respect to the compensated processes $\overline{X}$, $\overline{X}^{(k)}$, $\overline{\Phi}_{j}$ and $\overline{\Psi}_{i}^{(l)}$.
In fact, the above representation could also be derived for non-compensated processes at the cost of an additional Lebesgue integral with respect to an appropriate sum of compensators.
\end{remark}

\text{ } \text{ } This generalizes the classical results of Emery \cite{emery1989}, Dellacherie et al. \cite[p. 207]{Dell1992} and
Nualart and Schoutens \cite{NS}.\\

\section{Proof of predictable representation}\label{pol3}
\subsection{Polynomial representation}
The main result of this section gives a stochastic integral representation of polynomials
of the form $\overline{X}^{g}\cdot\overline{\Phi}^{p}_{j}\cdot\overline{\Psi}_{i}^{b}$.
We follow the argument of Nualart and Schoutens \cite{NS}.
\begin{theorem}\label{20}
We have the following representation:
\begin{eqnarray}\label{22}
 \nonumber \lefteqn{\overline{X}^{g}(t)\overline{\Phi}^{p}_{j}(t)\overline{\Psi}_{i}^{b}(t) =  f^{(g+p+b)}(t) +\sum_{s=1}^{g}\sum_{\tau=1}^{b}\sum_{\zeta=1}^{p}
 \sum_{\substack{(v_{1},\ldots,v_{\tau})\\  \in\{1,\ldots,b\}^{\tau}}}\sum_{\substack{(\iota_{1},\ldots,\iota_{s})\\  \in\{1,\ldots,g\}^{s}}}
  \int_{0}^{t} \int_{0}^{t_{1}-} \cdots \int_{0}^{t_{s+\tau+\zeta-1-}}}\\
 & & f^{(g+p+b)}_{(v_1,\dots, v_\tau, \iota_{1},\ldots,\iota_{s},i,j)}(t, t_{1},t_{2},\ldots,t_{s+\tau+\zeta})\mathrm{d}\overline{\Phi}_{j}(t_{s+\tau+\zeta})\ldots \mathrm{d}\overline{\Phi}_{j}(t_{s+\tau+1})\\
\nonumber & &\quad \quad \mathrm{d}\overline{\Psi}^{(v_\tau)}_{i}(t_{s+\tau})\ldots \mathrm{d}\overline{\Psi}^{(v_2)}_{i}(t_{s+2}) \mathrm{d}\overline{\Psi}^{(v_1)}_{i}(t_{s+1})\mathrm{d}\overline{X}^{(\iota_{s})}(t_{s}) \ldots \mathrm{d}\overline{X}^{(\iota_{2})}(t_{2})\mathrm{d}\overline{X}^{(\iota_{1})}(t_{1}),
\end{eqnarray}
where $f^{(g+p+b)}$ and $f^{(g+p+b)}_{(v_1,\ldots, v_\tau, \iota_{1},\ldots,\iota_{s},i,j)}$ are some random fields being a sum of products of predictable processes with respect to $t, t_{1},t_{2},\ldots,t_{s+\tau+\zeta}$   defined on $L^{2}(\Omega,\mathcal{F})$. 
\end{theorem}

\proof
We will express $\overline{X}^{g}(t)\overline{\Phi}^{p}_{j}(t)\overline{\Psi}^{b}_i(t)$ (for $t\geq 0$, $g, p, b \geq 0$, $i,j=1,\ldots,N$)
as a sum of stochastic integrals of lower powers of $\overline{X}, \overline{\Phi}_{j} $ and $\overline{\Psi}_{i}$
 with respect to the processes $\overline{X}^{(k)}$, $\overline{\Phi}_{j}$ and $\overline{\Psi}_{i}^{(l)}$ (for $k\leq g$ and $l\leq b$).
Note that the processes $\overline{\Phi}_{j}$ and $\overline{\Psi}_{i}$ have bounded variation and they are constant between jumps. Then (see Protter \cite[p. 75]{P}) the following holds true: $[\overline{X},\overline{\Phi}_{j}]^{c}(s)=0$, $[\overline{\Phi}_{j},\overline{\Phi}_{j}]^{c}(s)=0$, $[\overline{\Psi}_{i},\overline{\Phi}_{j}]^{c}(s)=0$, $[\overline{\Psi}_{i}, \overline{\Psi}_{i}]^{c}(s)=0$ and $[\overline{\Psi}_{i},\overline{X}]^{c}(s)=0$. Moreover, from the definition of $\overline{X}$ in  (\ref{1010}) we have  $[\overline{X},\overline{X}]^{c}(s)=\int_{0}^{s}\sigma_{0}^{2}(u-)\mathrm{d}u$.
Using It\^{o}'s formula (see Protter \cite[p. 81]{P}) we can write
\begin{eqnarray}\label{13}
\lefteqn{\overline{X}^{g}(t)\overline{\Phi}^{p}_{j}(t)\overline{\Psi}_{i}^{b}(t) = \overline{X}^{g}(0)\overline{\Phi}^{p}_{j}(0)\overline{\Psi}_{i}^{b}(0)+\int_{0}^{t}g\overline{X}^{g-1}(s-)\overline{\Phi}^{p}_{j}(s-)\overline{\Psi}_{i}^{b}(s-)\mathrm{d}\overline{X}(s)}\\
\nonumber & & + \int_{0}^{t}p\overline{X}^{g}(s-)\overline{\Phi}^{p-1}_{j}(s-)\overline{\Psi}_{i}^{b}(s-)\mathrm{d}\overline{\Phi}_{j}(s) +\int_{0}^{t}b\overline{X}^{g}(s-)\overline{\Phi}^{p}_{j}(s-)\overline{\Psi}_{i}^{b-1}(s-)\mathrm{d}\overline{\Psi}_{i}(s) \\
\nonumber & & + \frac{1}{2}g(g-1) \mathbf{I}_{1}+\mathbf{I}_{2},
\end{eqnarray}
where 
\begin{equation*}
\mathbf{I}_{1} := \int_{0}^{t}\overline{X}^{g-2}(s-)\overline{\Phi}^{p}_{j}(s-)\overline{\Psi}_{i}^{b}(s-)\sigma_{0}^{2}(s-)\mathrm{d}s
\end{equation*}
and
\begin{eqnarray}\label{i2}
\lefteqn{\mathbf{I}_{2}:=\sum_{0<s\leq t}\Big\{ \mathbf{I}_{3}
- \overline{X}^{g}(s-)\overline{\Phi}^{p}_{j}(s-)\overline{\Psi}_{i}^{b}(s-)-g\overline{X}^{g-1}(s-)\overline{\Phi}^{p}_{j}(s-)\overline{\Psi}_{i}^{b}(s-)\Delta \overline{X}(s)}\\
\nonumber & & - p\overline{X}^{g}(s-)\overline{\Phi}^{p-1}_{j}(s-)\overline{\Psi}_{i}^{b}(s-)\Delta \overline{\Phi}_{j}(s) -b  \overline{X}^{g}(s-)\overline{\Phi}^{p}_{j}(s-)\big(\overline{\Psi}_{i}\big)^{b-1}(s-)\Delta \overline{\Psi}_{i}(s)\Big\}
\end{eqnarray}
with
\begin{equation}\label{i3}
\mathbf{I}_{3} := \overline{X}^{g}(s)\overline{\Phi}^{p}_{j}(s)\overline{\Psi}_{i}^{b}(s).
\end{equation}
Using integration by parts we can rewrite $\mathbf{I}_{1}$ as follows:
\begin{eqnarray}\label{12} 
\lefteqn{ \mathbf{I}_{1}=  \overline{X}^{g-2}(t-)\overline{\Phi}^{p}_{j}(t-)\overline{\Psi}_{i}^{b}(t-)\int_{0}^{t}\sigma_{0}^{2}(u-)\mathrm{d}u -\int_{0}^{t}\Bigg(\int_{0}^{s}\sigma_{0}^{2}(u-)\mathrm{d}u\Bigg)\overline{X}^{g-2}(s)\overline{\Phi}^{p}_{j}(s)\mathrm{d}\overline{\Psi}_{i}^{b}(s)}\nonumber\\
& &  - \int_{0}^{t}\Bigg(\int_{0}^{s}\sigma_{0}^{2}(u-)\mathrm{d}u\Bigg) \overline{X}^{g-2}(s)\overline{\Psi}_{i}^{b}(s)\mathrm{d}\overline{\Phi}_{j}(s)-\int_{0}^{t}\Bigg(\int_{0}^{s}\sigma_{0}^{2}(u-)\mathrm{d}u\Bigg) \overline{X}^{g-2}(s)\mathrm{d}[\overline{\Phi}^{p}_{j},\overline{\Psi}_{i}^{b}](s)\\
\nonumber & & - \int_{0}^{t}\Bigg(\int_{0}^{s}\sigma_{0}^{2}(u-)\mathrm{d}u\Bigg) \overline{\Phi}^{p}_{j}(s)\overline{\Psi}_{i}^{b}(s)\mathrm{d}\overline{X}^{g-2}(s)
\end{eqnarray}
with
\[\mathrm{d}[\overline{\Phi}^{p}_{j},\overline{\Psi}_{i}^{b}](s)=\mathrm{d}\overline{\Psi}_{i}^{b}(s)\quad\text{if $i=j$} \]
and $[\overline{\Phi}^{p}_{j},\overline{\Psi}_{i}^{b}](s)=0$ otherwise.
Note that $\overline{X}(s)=\overline{X}(s-)+\Delta \overline{X}(s)$, $\overline{\Phi}_{j}(s)=\overline{\Phi}_{j}(s-)+\Delta \overline{\Phi}_{j}(s)$ and $\overline{\Psi}_{i}(s)=\overline{\Psi}_{i}(s-)+\Delta \overline{\Psi}_{i}(s)$.
Using the Binomial Theorem we can rewrite $\mathbf{I}_{3} $, defined in (\ref{i3}), as follows: 
\begin{eqnarray}
\nonumber \mathbf{I}_{3} & = & \big(\overline{X}(s-)+\Delta \overline{X}(s)\big)^{g}\big(\overline{\Phi}_{j}(s-)+\Delta \overline{\Phi}_{j}(s)\big)^{p} \big(\overline{\Psi}_{i}(s-)+\Delta \overline{\Psi}_{i}(s)\big)^{b}\\
 \nonumber & = & \mathbf{I}_{4} \cdot \mathbf{I}_{5} \cdot  \mathbf{I}_{6},
\end{eqnarray}
where
\begin{equation}\label{i4}
\mathbf{I}_{4} := \overline{X}^{g}(s-)+\sum_{m_{1}=1}^{g} {g\choose m_{1}}\overline{X}^{g-m_{1}}(s-)\big(\Delta \overline{X}(s) \big)^{m_{1}},
\end{equation}
\begin{equation}\label{i5}
\mathbf{I}_{5} :=\overline{\Phi}_{j}^{p}(s-)+\sum_{m_{2}=1}^{p} {p\choose m_{2}} \overline{\Phi}_{j}^{p-m_{2}}(s-)\big(\Delta \overline{\Phi}_{j}(s) \big)^{m_{2}},
\end{equation}
\begin{equation}\label{i6}
\mathbf{I}_{6} :=\overline{\Psi}_{i}^{b}(s-)+ \sum_{m_{3}=1}^{b} {b\choose m_{3}}\overline{\Psi}_{i}^{b-m_{3}}(s-)\big(\Delta \overline{\Psi}_{i}(s) \big)^{m_{3}}.
\end{equation}
Then
\begin{eqnarray}\label{i4i61}
\nonumber \lefteqn{  \mathbf{I}_{4}  \cdot \mathbf{I}_{6} =  \Bigg(\sum_{m_{1}=1}^{g} {g\choose m_{1}}\overline{X}^{g-m_{1}}(s-)\big(\Delta \overline{X}(s) \big)^{m_{1}}\Bigg) \Bigg(\sum_{m_{3}=1}^{b} {b\choose m_{3}}\overline{\Psi}_{i}^{b-m_{3}}(s-)\big(\Delta \overline{\Psi}_{i}(s) \big)^{m_{3}}\Bigg)+ \overline{X}^{g}(s-)\overline{\Psi}_{i}^{b}(s-)}\\
& & + \sum_{m_{1}=1}^{g} {g\choose m_{1}}\overline{X}^{g-m_{1}}(s-)\overline{\Psi}_{i}^{b}(s-)\big(\Delta \overline{X}(s) \big)^{m_{1}}  + \sum_{m_{3}=1}^{b} {b\choose m_{3}} \overline{X}^{g}(s-)\overline{\Psi}_{i}^{b-m_{3}}(s-)\big(\Delta \overline{\Psi}_{i}(s) \big)^{m_{3}}.
\end{eqnarray}
Since $\overline{X}$ and $\overline{\Psi}_{i}$ do not jump at the same time, we have
 $\Delta \overline{X}(s)\Delta \overline{\Psi}_{i}(s)=0 $. Thus, the first component of the right-hand side of (\ref{i4i61}) is zero, and so
\begin{eqnarray}\label{i4i6}
\lefteqn{   \mathbf{I}_{4}  \cdot \mathbf{I}_{6}  =  \overline{X}^{g}(s-)\overline{\Psi}_{i}^{b}(s-)+\sum_{m_{1}=1}^{g} {g\choose m_{1}}\overline{X}^{g-m_{1}}(s-)\overline{\Psi}_{i}^{b}(s-)\big(\Delta \overline{X}(s) \big)^{m_{1}}  }\\
\nonumber & & + \sum_{m_{3}=1}^{b} {b\choose m_{3}} \overline{X}^{g}(s-)\overline{\Psi}_{i}^{b-m_{3}}(s-)\big(\Delta \overline{\Psi}_{i}(s) \big)^{m_{3}}.
\end{eqnarray}
By (\ref{i5}) and (\ref{i4i6}),
\begin{eqnarray}\label{i4i5i6}
\lefteqn{\mathbf{I}_{3}= \mathbf{I}_{4} \cdot \mathbf{I}_{5} \cdot \mathbf{I}_{6}  = \overline{X}^{g}(s-)\overline{\Phi}^p_{j}(s-)\overline{\Psi}_{i}^{b}(s-)+ \sum_{m_{2}=1}^{p} {p\choose m_{2}} \overline{X}^{g}(s-) \overline{\Phi}_{j}^{p-m_{2}}(s-)\overline{\Psi}_{i}^{b}(s-)\big(\Delta \overline{\Phi}_{j}(s) \big)^{m_{2}}}\nonumber \\
& & +\Bigg(\sum_{m_{1}=1}^{g} {g\choose m_{1}}\overline{X}^{g-m_{1}}(s-)\overline{\Psi}_{i}^{b}(s-)\big(\Delta \overline{X}(s) \big)^{m_{1}}\Bigg)  \Bigg(\sum_{m_{2}=1}^{p} {p\choose m_{2}} \overline{\Phi}_{j}^{p-m_{2}}(s-)\big(\Delta \overline{\Phi}_{j}(s) \big)^{m_{2}}\Bigg)\\
\nonumber & & + \Bigg(\sum_{m_{3}=1}^{b} {b\choose m_{3}}\overline{X}^{g}(s-)\overline{\Psi}_{i}^{b-m_{3}}(s-)\big(\Delta \overline{\Psi}_{i}(s) \big)^{m_{3}}\Bigg) \Bigg(\sum_{m_{2}=1}^{p} {p\choose m_{2}} \overline{\Phi}_{j}^{p-m_{2}}(s-)\big(\Delta \overline{\Phi}_{j}(s) \big)^{m_{2}}\Bigg)\\
\nonumber & & +\Bigg(\sum_{m_{1}=1}^{g} {g\choose m_{1}}\overline{X}^{g-m_{1}}(s-)\overline{\Phi}^p_{j}(s-)\overline{\Psi}_{i}^{b}(s-)\big(\Delta \overline{X}(s) \big)^{m_{1}}\Bigg)\\
\nonumber & & + \Bigg(\sum_{m_{3}=1}^{b} {b\choose m_{3}}\overline{X}^{g}(s-)\overline{\Phi}^p_{j}(s-)\overline{\Psi}_{i}^{b-m_{3}}(s-)\big(\Delta \overline{\Psi}_{i}(s) \big)^{m_{3}}\Bigg).
\end{eqnarray}
The third component of the above sum is zero. This follows from the observation that $\Delta \overline{X}(s)\Delta \overline{\Phi}_{j}(s)=0$, because $\overline{X}$ and $\overline{\Phi}_{j}$ do not jump at the same time. Note that the fifth component of the sum (\ref{i4i5i6}) can be written as follows:
\begin{eqnarray}\label{Bigg}
\lefteqn{\Bigg(\sum_{m_{1}=1}^{g} {g\choose m_{1}}\overline{X}^{g-m_{1}}(s-)\overline{\Phi}^p_{j}(s-)\overline{\Psi}_{i}^{b}(s-)\big(\Delta \overline{X}(s) \big)^{m_{1}}\Bigg)}\nonumber \\
& &= g\overline{X}^{g-1}(s-)\overline{\Phi}^p_{j}(s-)\overline{\Psi}_{i}^{b}(s-) \Delta \overline{X}(s) +\sum_{m_{1}=2}^{g} {g\choose m_{1}}\overline{X}^{g-m_{1}}(s-)\overline{\Phi}^p_{j}(s-)\overline{\Psi}_{i}^{b}(s-)\big(\Delta \overline{X}(s) \big)^{m_{1}}.
\end{eqnarray}
Combining (\ref{i4i5i6}) and (\ref{Bigg}) we conclude that 
\begin{eqnarray}\label{i4i6bdds}
\lefteqn{\mathbf{I}_{3}= \overline{X}^{g}(s-)\overline{\Phi}^p_{j}(s-)\overline{\Psi}_{i}^{b}(s-)+  \sum_{m_{2}=1}^{p} {p\choose m_{2}} \overline{X}^{g}(s-) \overline{\Phi}_{j}^{p-m_{2}}(s-)\overline{\Psi}_{i}^{b}(s-)\big(\Delta \overline{\Phi}_{j}(s) \big)^{m_{2}}}\nonumber \\
\nonumber & & + \Bigg(\sum_{m_{3}=1}^{b} {b\choose m_{3}}\overline{X}^{g}(s-)\overline{\Psi}_{i}^{b-m_{3}}(s-)\big(\Delta \overline{\Psi}_{i}(s) \big)^{m_{3}}\Bigg) \Bigg(\sum_{m_{2}=1}^{p} {p\choose m_{2}} \overline{\Phi}_{j}^{p-m_{2}}(s-)\big(\Delta \overline{\Phi}_{j}(s) \big)^{m_{2}}\Bigg)\\
\nonumber & &+g\overline{X}^{g-1}(s-)\overline{\Phi}^p_{j}(s-)\overline{\Psi}_{i}^{b}(s-) \Delta \overline{X}(s) +\sum_{m_{1}=2}^{g} {g\choose m_{1}}\overline{X}^{g-m_{1}}(s-)\overline{\Phi}^p_{j}(s-)\overline{\Psi}_{i}^{b}(s-)\big(\Delta \overline{X}(s) \big)^{m_{1}}\\
\nonumber & & + \sum_{m_{3}=1}^{b} {b\choose m_{3}}\overline{X}^{g}(s-)\overline{\Phi}^p_{j}(s-)\overline{\Psi}_{i}^{b-m_{3}}(s-)\big(\Delta \overline{\Psi}_{i}(s) \big)^{m_{3}}.
\end{eqnarray}
 Now, inserting the above $\mathbf{I}_{3}$ into $\mathbf{I}_{2}$ defined by (\ref{i2}), we derive
\begin{eqnarray}
\lefteqn{ \mathbf{I}_{2} =  \sum_{0<s\leq t}\Bigg\{   \sum_{m_{2}=1}^{p} {p\choose m_{2}} \overline{X}^{g}(s-) \overline{\Phi}_{j}^{p-m_{2}}(s-)\overline{\Psi}_{i}^{b}(s-)\big(\Delta \overline{\Phi}_{j}(s) \big)^{m_{2}}}\nonumber \\
 \nonumber & & + \Bigg(\sum_{m_{3}=1}^{b} {b\choose m_{3}}\overline{X}^{g}(s-)\overline{\Psi}_{i}^{b-m_{3}}(s-)\big(\Delta \overline{\Psi}_{i}(s) \big)^{m_{3}}\Bigg) \Bigg(\sum_{m_{2}=1}^{p} {p\choose m_{2}} \overline{\Phi}_{j}^{p-m_{2}}(s-)\big(\Delta \overline{\Phi}_{j}(s) \big)^{m_{2}}\Bigg)\\
\nonumber & & + \sum_{m_{1}=2}^{g} {g\choose m_{1}} \overline{X}^{g-m_{1}}(s-)\overline{\Phi}_{j}^{p}(s-) \overline{\Psi}_{i}^{b}(s-)\big(\Delta \overline{X}(s) \big)^{m_{1}} + \sum_{m_{3}=1}^{b} {b\choose m_{3}}\overline{X}^{g}(s-)\overline{\Phi}^{p}_{j}(s-)\overline{\Psi}_{i}^{b-m_{3}}(s-)\big(\Delta \overline{\Psi}_{i}(s) \big)^{m_{3}}\\
\nonumber & & - p\overline{X}^{g}(s-)\overline{\Phi}^{p-1}_{j}(s-)\overline{\Psi}_{i}^{b}(s-)\Delta \overline{\Phi}_{j}(s) -b  \overline{X}^{g}(s-)\overline{\Phi}^{p}_{j}(s-)\big(\overline{\Psi}_{i}\big)^{b-1}(s-)\Delta \overline{\Psi}_{i}(s)\Bigg\}.
\end{eqnarray}
Note that $\big(\Delta \overline{\Phi}_{j}(s) \big)^{m_{2}}= \Delta \overline{\Phi}_{j}(s)= \Delta \Phi_{j}(s) $ and $\Delta \overline{\Psi}_{i}(s) = \Delta \Psi_{i}(s) $, which follows from the definition of $\overline{\Phi}_{j}$ and $\overline{\Psi}_{i}$ in (\ref{ovlPhi}) and (\ref{ovlPsi}), respectively. Moreover, $ \Delta \Phi_{j}(s) \Delta \Psi_{i}(s)=\delta_{ij}  \Delta \Psi_{i}(s)$ since $\Phi_{j}$ and $\Psi_{i}$  jump at the same time only when $i=j$. If $i \neq j$ then either $ \Delta \Psi_{i}$ or $\Delta \Phi_{j}$ is zero. Thus we can rewrite $\mathbf{I}_{2}$ as follows:
\begin{eqnarray*}
\lefteqn{ \nonumber \mathbf{I}_{2} =  \sum_{0<s\leq t}\Bigg\{  \sum_{m_{2}=1}^{p} {p\choose m_{2}} \overline{X}^{g}(s-)  \overline{\Phi}_{j}^{p-m_{2}}(s-)\overline{\Psi}_{i}^{b}(s-) \Delta \Phi_{j}(s)}\\
 & & + \Bigg(\sum_{m_{3}=1}^{b} {b\choose m_{3}}\overline{X}^{g}(s-)\overline{\Psi}_{i}^{b-m_{3}}(s-)\Bigg) \Bigg(\sum_{m_{2}=1}^{p} {p\choose m_{2}} \overline{\Phi}_{j}^{p-m_{2}}(s-)\Bigg)\delta_{ij}  \big(\Delta \Psi_{i}(s) \big)^{m_{3}}\\
\nonumber & &+ \sum_{m_{1}=2}^{g} {g\choose m_{1}} \overline{X}^{g-m_{1}}(s-) \overline{\Phi}^{p}_{j}(s-)\overline{\Psi}_{i}^{b}(s-)\big(\Delta \overline{X}(s) \big)^{m_{1}}+ \sum_{m_{3}=1}^{b} {b\choose m_{3}}\overline{X}^{g}(s-)\overline{\Psi}_{i}^{b-m_{3}}(s-)\big(\Delta \Psi_{i}(s) \big)^{m_{3}}\\
\nonumber & & - p\overline{X}^{g}(s-)\overline{\Phi}^{p-1}_{j}(s-)\overline{\Psi}_{i}^{b}(s-)\Delta \Phi_{j}(s) -b  \overline{X}^{g}(s-)\overline{\Phi}^{p}_{j}(s-)\big(\overline{\Psi}_{i}\big)^{b-1}(s-)\Delta \Psi_{i}(s)\Bigg\}.
\end{eqnarray*}
Now, we rewrite $\mathbf{I}_{2}$ as a sum of some stochastic integrals: 
\begin{eqnarray*}
\lefteqn{ \mathbf{I}_{2} =  \sum_{m_{2}=1}^{p}\int_0^t {p\choose m_{2}} \overline{X}^{g}(s-)\overline{\Phi}_{j}^{p-m_{2}}(s-)\overline{\Psi}_{i}^{b}(s-)  \mathrm{d} \Phi_{j}(s)}\nonumber\\
\nonumber & & + \Bigg(\sum_{m_{3}=1}^{b} \int_0^t {b\choose m_{3}}\overline{X}^{g}(s-) \overline{\Psi}_{i}^{b-m_{3}}(s-)\Bigg)\Bigg(\sum_{m_{2}=1}^{p} {p\choose m_{2}} \overline{\Phi}_{j}^{p-m_{2}}(s-) \Bigg)\delta_{ij} \mathrm{d} \Psi_{i}^{(m_{3})}(s)\\
\nonumber  & & + \sum_{m_{1}=2}^{g} \int_0^t {g\choose m_{1}} \overline{X}^{g-m_{1}}(s-) \overline{\Phi}^{p}_{j}(s-)\overline{\Psi}_{i}^{b}(s-) \mathrm{d} X^{(m_{1})}(s) +\sum_{m_{3}=1}^{b} \int_0^t {b\choose m_{3}}\overline{X}^{g}(s-)\overline{\Psi}_{i}^{b-m_{3}}(s-) \mathrm{d} \Psi_{i}^{(m_{3})} (s)\\
\nonumber & & - \int_0^t p\overline{X}^{g}(s-)\overline{\Phi}^{p-1}_{j}(s-)\overline{\Psi}_{i}^{b}(s-)\mathrm{d} \Phi_{j}(s)  - \int_0^t b\overline{X}^{g}(s-)\overline{\Phi}^{p}_{j}(s-)\overline{\Psi}_{i}^{b-1}(s-) \mathrm{d} \Psi_{i}(s).
\end{eqnarray*}
We shall rewrite the above expression in terms of integrals with respect to the compensated processes $\overline{X}^{(m_{1})}$,  $\overline{\Phi}_{j}$ and $\overline{\Psi}_{i}$ given in (\ref{ovlX}), (\ref{ovlPhi}) and (\ref{ovlPsi}), respectively:
\begin{eqnarray}\label{1444}
\lefteqn{  \mathbf{I}_{2} =  \sum_{m_{2}=1}^{p} \int_0^t  {p\choose m_{2}} \overline{X}^{g}(s-)\overline{\Phi}_{j}^{p-m_{2}}(s-)\overline{\Psi}_{i}^{b}(s-)  \mathrm{d} \overline{\Phi}_{j}(s)}\nonumber\\
& & +\sum_{m_{3}=1}^{b} \int_0^t  \Bigg( {b\choose m_{3}}\overline{X}^{g}(s-) \overline{\Psi}_{i}^{b-m_{3}}(s-)\Bigg)\Bigg(\sum_{m_{2}=1}^{p} {p\choose m_{2}} \overline{\Phi}_{j}^{p-m_{2}}(s-) \Bigg)\delta_{ij} \mathrm{d} \overline{\Psi}_{i}^{(m_{3})}(s)\\
\nonumber & &+ \sum_{m_{1}=2}^{g} \int_0^t  {g\choose m_{1}} \overline{X}^{g-m_{1}}(s-) \overline{\Phi}^{p}_{j}(s-)\overline{\Psi}_{i}^{b}(s-) \mathrm{d} \overline{X}^{(m_{1})}(s)\\
\nonumber & & +\sum_{m_{3}=1}^{b} \int_0^t {b\choose m_{3}}\overline{X}^{g}(s-)\overline{\Psi}_{i}^{b-m_{3}}(s-) \mathrm{d} \overline{\Psi}_{i}^{(m_{3})} (s) \\
\nonumber & & - \int_0^t p\overline{X}^{g}(s-)\overline{\Phi}^{p-1}_{j}(s-)\overline{\Psi}_{i}^{b}(s-)\mathrm{d} \overline{\Phi}_{j}(s)- \int_0^t b\overline{X}^{g}(s-)\overline{\Phi}^{p}_{j}(s-)\overline{\Psi}_{i}^{b-1}(s-)\mathrm{d} \overline{\Psi}_{i}(s)+ \mathbf{I}_{7} ,
\end{eqnarray}
where
\begin{eqnarray}\label{1441}
\lefteqn{   \mathbf{I}_{7} := \sum_{m_{2}=1}^{p} \int_0^t  {p\choose m_{2}} \overline{X}^{g}(s-)\overline{\Phi}_{j}^{p-m_{2}}(s-)\overline{\Psi}_{i}^{b}(s-)  \lambda_{j}(s)\mathrm{d}s}\nonumber\\ 
& & +\sum_{m_{3}=1}^{b} \int_0^t \Bigg( {b\choose m_{3}}\overline{X}^{g}(s-) \overline{\Psi}_{i}^{b-m_{3}}(s-)\Bigg)\Bigg(\sum_{m_{2}=1}^{p} {p\choose m_{2}} \overline{\Phi}_{j}^{p-m_{2}}(s-) \Bigg)\delta_{ij} \mathbb{E}\big(U^{(j)}_{n}\big)^{m_3}\lambda_{j}(s)\mathrm{d}s\\
  & & +   \sum_{m_{1}=2}^{g} \int_0^t \int\limits_{\mathbb{R}} {g\choose m_{1}} \overline{X}^{g-m_{1}}(s-)\overline{\Phi}^{p}_{j}(s-)\overline{\Psi}_{i}^{b}(s-) \gamma^{m_1}(s-,x)\nu(\mathrm{d}x) \mathrm{d}s\nonumber\\
  \nonumber & & +\sum_{m_{3}=1}^{b} \int_0^t {b\choose m_{3}}\overline{X}^{g}(s-)\overline{\Psi}_{i}^{b-m_{3}}(s-) \delta_{ij} \mathbb{E}\big(U^{(j)}_{n}\big)^{m_3}\lambda_{j}(s)\mathrm{d}s\\
 \nonumber & & - \int_0^t p\overline{X}^{g}(s-)\overline{\Phi}^{p-1}_{j}(s-)\overline{\Psi}_{i}^{b}(s-)\lambda_{j}(s)\mathrm{d}s- \int_0^t b\overline{X}^{g}(s-)\overline{\Phi}^{p}_{j}(s-)\overline{\Psi}_{i}^{b-1}(s-)\delta_{ij} \mathbb{E}\big(U^{(i)}_{n}\big)\lambda_{i}(s)\mathrm{d}s.
\end{eqnarray} 
Note that $\mathbf{I}_{7}$ is a sum of terms of the form
\begin{equation}\label{i8}
\mathbf{I}_{8} :=  \int_0^t  \overline{X}^{g-m_1}(s-)\overline{\Phi}^{p-m_2}_{j}(s-)\overline{\Psi}_{i}^{b-m_3}(s-) \Theta(s-)\mathrm{d}s
\end{equation}
for various  predictable processes $ \Theta(s-)$ (i.e. $ \Theta(s-)= \lambda_{j}(s)$ or $ \Theta(s-)=\delta_{ij} \mathbb{E}\big(U^{(j)}_{n}\big)^{m_3}\lambda_{j}(s)$ for $m_3=1,...,b$ or $ \Theta(s-)= \int_{\mathbb{R}} \gamma^{m_1}(s-,x)\nu(\mathrm{d}x)$ for $m_1=2,...,g$ ).
Similarly, as in (\ref{12}) using integration by parts for semimartingales we can rewrite (\ref{i8}) in the following way:
\begin{eqnarray}\label{integrbypart}
\lefteqn{ \mathbf{I}_{8}=\overline{X}^{g-m_1}(t-)\overline{\Phi}^{p-m_2}_{j}(t-)\overline{\Psi}_{i}^{b-m_3}(t-)\int_{0}^{t} \Theta(u-)\mathrm{d}u  }\nonumber\\
& & -\int_{0}^{t}\Bigg(\int_{0}^{s} \Theta(u-)\mathrm{d}u\Bigg)\overline{X}^{g-m_1}(s)\overline{\Phi}^{p-m_2}_{j}(s)\mathrm{d}\overline{\Psi}_{i}^{b-m_3}(s)\\
\nonumber  & &  - \int_{0}^{t}\Bigg(\int_{0}^{s} \Theta(u-)\mathrm{d}u\Bigg) \overline{X}^{g-m_1}(s)\overline{\Psi}_{i}^{b-m_3}(s)\mathrm{d}\overline{\Phi}_{j}(s)\\
\nonumber & & -\int_{0}^{t}\Bigg(\int_{0}^{s} \Theta(u-)\mathrm{d}u\Bigg) \overline{X}^{g-m_1}(s)\mathrm{d}[\overline{\Phi}^{p-m_2}_{j},\overline{\Psi}_{i}^{b-m_3}](s)\\
\nonumber & & - \int_{0}^{t}\Bigg(\int_{0}^{s} \Theta(u-)\mathrm{d}u\Bigg) \overline{\Phi}^{p-m_2}_{j}(s)\overline{\Psi}_{i}^{b-m_3}(s)\mathrm{d}\overline{X}^{g-m_1}(s).
\end{eqnarray}
Combining (\ref{12}), (\ref{1444}), (\ref{1441}) and (\ref{integrbypart})  we can express (\ref{13}) as a sum of stochastic integrals of powers strictly lower than $g+p+b$ of $\overline{X}, \overline{\Phi}_{j} $ and $\overline{\Psi}_{i}$
 with respect to the processes $\overline{X}^{(k)}$, $\overline{\Phi}_{j}$ and $\overline{\Psi}_{i}^{(l)}$ for $k\leq g$ and $l\leq b$.
Then all the above identities and induction complete the proof.\\
\halmos

\section{Proofs of main results}\label{sec:proos}
\subsection{Proof of Theorem \ref{987}}\label{pol1}
We start with the following proposition.

\begin{proposition}\label{apro}
For fixed $t\geq 0$ and $s_1\leq \ldots \leq s_m\leq t$ let
\begin{eqnarray*}
\boldsymbol{\overline{X}}&:=&(\overline{X}(s_{1}), \ldots , \overline{X}(s_{m}) ),\\
\boldsymbol{\overline{\Phi}}_{j}&:=&(\overline{\Phi}_{j}(s_{1}), \ldots , \overline{\Phi}_{j}(s_{m}) ) \quad\text{ for } j=1,...,N,\\
\boldsymbol{\overline{\Psi}}_{i} &:=&(\overline{\Psi}_{i}(s_{1}), \ldots , \overline{\Psi}_{i}(s_{m})) \quad\text{ for } i=1,...,N
\end{eqnarray*}
and
\begin{equation*}
Z\in L^{2}\big(\Omega, \mathcal{F}_t \big).
\end{equation*}
Then, for any $\varepsilon>0 $, there exists $m \in  \mathbb{N} $ and a random variable $Z_{\varepsilon}\in L^{2}\Big(\Omega, \sigma \big( \boldsymbol{\overline{X}}, \boldsymbol{\overline{\Phi}}_{1},\ldots, \boldsymbol{\overline{\Phi}}_{N}, \boldsymbol{\overline{\Psi}}_{1}, \ldots, \boldsymbol{\overline{\Psi}}_{N}\big)\Big)$
such that
\begin{equation*}
\mathbb{E}\big[(Z-Z_{\varepsilon})^{2} \big]<\varepsilon.
\end{equation*}
\end{proposition}
\proof
The conclusion follows from more general considerations.
Let $Z\in L^2(\Omega,\sigma(Z_1,Z_2,\ldots))$ for certain variables $\{Z_1,Z_2,\dots\}$.
Then $Z$ can be approximated with any given accuracy in $L^2$ norm by $$\tilde{Z}_\epsilon:=\sum\limits_{l=1}^{m} k_l \textbf{1}_{A_l},$$
where $A_l=\left\{Z_1\in \mathcal{B}_1^l, Z_2\in \mathcal{B}_2^l,\ldots\right\}$ for Borel sets $\mathcal{B}_k^l$ and constants $k_l$.\\
For any $\epsilon>0$, there exists a positive integer $M$ such that $$ \mathbb{P}(A_l^M\setminus A_l)=\mathbb{P}(Z_{M+1}\in \mathcal{B}_{M+1}^l, Z_{M+2}\in \mathcal{B}_{M+2}^l, \ldots)=\mathbb{P}(\bigcap \limits_{n=M+1}^{\infty}Z_{n}\in \mathcal{B}_{n}^l)\leq \epsilon,$$
 where
$A_l^M=\left\{Z_1\in \mathcal{B}_1^l, \ldots, Z_M\in \mathcal{B}_M^l \right\}$. 
It remains to define $Z_\epsilon:= \sum\limits_{l=1}^{m}k_l \textbf{1}_{A_l^M}$.\\
\halmos\\
Let
\begin{eqnarray*}
\mathcal{P} &:=&  \Big\{\overline{X}^{g_{1}}(t_{1})\cdot \ldots \cdot\overline{X}^{g_{m}} (t_{m}) \cdot \prod_{i,j=1}^{N} \overline{\Phi}^{p_{1}}_{j}(t_{1})\cdot \ldots \cdot  \overline{\Phi}_{j}^{p_{m}}(t_{m}) \overline{\Psi}_{i}^{b_{1}}(t_{1})  \cdot \ldots  \cdot \overline{\Psi}_{i}^{b_{m}}(t_{m}):0\leq t_{1}< \ldots\\
\nonumber & & < t_{m}, g_{1}, \ldots ,g_{m} \geq 0, p_{1}, \ldots ,p_{m} \geq 0, b_{1}, \ldots ,b_{m}\geq 0\Big\}.
\end{eqnarray*}
\begin{lemma}\label{totalfamily}
The set $\mathcal{P}$ is a total family in $L^{2}(\Omega, \mathcal{F}) $, i.e. the linear span of $\mathcal{P}$ is dense in $L^{2}(\Omega, \mathcal{F}) $.
\end{lemma}
\proof
Assume that $\mathcal{P}$ is not a total family. Then there is $Z\in L^{2}(\Omega, \mathcal{F}) $ such that $Z$ is orthogonal to $\mathcal{P}$.
From Proposition \ref{apro} there exists a Borel function $f$ such that
\begin{equation*}
Z_{\varepsilon}=f_{\varepsilon}\big( \boldsymbol{\overline{X}},  \boldsymbol{\overline{\Phi}}_{1},\ldots,  \boldsymbol{\overline{\Phi}}_{N},  \boldsymbol{\overline{\Psi}}_{1}, \ldots , \boldsymbol{\overline{\Psi}}_{N} \big).
\end{equation*}
\text{ } \text{ } Recall that under assumption \eqref{42} the polynomials are dense in $L^{2}(\mathbb{R}, \mathrm{d}\varphi)$, so we can approximate $Z_{\varepsilon}$ by polynomials. Furthermore, since
$Z\perp \mathcal{P}$, we have $\mathbb{E}(ZZ_{\varepsilon})=0$. Then from the Schwarz inequality we obtain
\begin{equation*}
\mathbb{E}\big(Z^{2}\big)=\mathbb{E}\big(Z(Z-Z_{\varepsilon})\big)\leq \sqrt{\mathbb{E}(Z^{2})\mathbb{E}\big((Z-Z_{\varepsilon})^{2}\big)}\leq \sqrt{\varepsilon \mathbb{E}(Z^{2})}.
\end{equation*}
Letting $\varepsilon\longrightarrow 0 $ yields $Z=0$ a.s.\\
\halmos\\

\text{ } \text{ } {\bf Proof of Theorem \ref{987}}.
Since by a linear transformation we can switch from the $\overline{X}^{(k)}$ to the $H^{(k)}$, from the $\overline{\Psi}_{i}^{(n)}$  to the $G_{i}^{(n)}$ and from the $\overline{\Phi}_{i}$ to the $G_{i}^{(1)}$, we can rewrite representation (\ref{22}) as follows:
\begin{eqnarray}\label{22bis}
 \nonumber \lefteqn{\overline{X}^{g}(t)\overline{\Phi}^{p}_{j}(t)\overline{\Psi}_{i}^{b}(t) =  f^{(g+p+b)}(t) +\sum_{s=1}^{g}\sum_{\tau=1}^{b+p}  \sum_{\iota_{1},\ldots,\iota_{s}\geq 1
} \sum_{\upsilon_{1},\ldots,\upsilon_{\tau}\geq 1} \int_{0}^{t} \int_{0}^{t_{1}-}}\\
& & \cdots \int_{0}^{t_{s+\tau-1-}} f^{(g+p+b)}_{(\upsilon_{1},\ldots,\upsilon_{\tau}, \iota_{1},\ldots,\iota_{s}, i)}(t_{1},t_{2},\ldots,t_{s+\tau })\\
 \nonumber & & \qquad\qquad \mathrm{d}G^{(\upsilon_{\tau})}_{i}(t_{s+\tau}) \ldots \mathrm{d}G^{(\upsilon_{1})}_{i}(t_{s+1})\mathrm{d}H^{(\iota_{s})}(t_{s}) \ldots \mathrm{d}H^{(\iota_{1})}(t_{1}),
\end{eqnarray}
where $ f^{(g+p+b)}$  and $f^{(g+p+b)}_{(\upsilon_{1},\ldots,\upsilon_{\tau}, \iota_{1},\ldots,\iota_{s}, i )}$ are random fields on $L^{2}(\Omega,\mathcal{F})$, and $G^{(l)}_{i}$ and $H^{(k)}$
$(i=1,\dots, N$,  $l,k\geq1)$ are orthogonal martingales defined in \eqref{defG} and \eqref{defH}, respectively.
In fact, stochastic integrals with respect to orthogonal martingales are again orthogonal (see Protter \cite[ Lemma 2, p. 180 and Theorem 35, p. 178]{P} ).\\
\text{ } \text{ } Since by Lemma \ref{totalfamily} the set $\mathcal{P}$ is a total family in $L^{2}(\Omega, \mathcal{F}) $, any $\mathcal{F}_t$-measurable
random variable $F$ in $L^{2}(\Omega,\mathcal{F}_{t})$ can be approximated by a polynomial constructed from the processes of the form $\overline{X}^{g}\cdot\overline{\Phi}^{p}_{j}\cdot\overline{\Psi}_{i}^{b}$
appearing on the left side of \eqref{22bis} where $i,j$ vary over $\{1,\ldots,N\}$. First observe that if $A$ and $B$ have the form of a multiple integral with respect to $G^{(l)}_{i}$ and $H^{(k)}$ $(i=1,\dots, N$,  $l,k\geq1)$ as on the right side of \eqref{22bis}, then $AB$ is also of this form (see the proof of \cite[Thm. 1]{NS} for details). Thus every random variable $F$ in $L^{2}(\Omega,\mathcal{F}_{t})$ has the following representation:
\begin{eqnarray*}
F(t) &=&  \mathbb{E}[F(t)]+\sum_{i=1}^{N}\sum_{s=1}^{\infty}\sum_{\tau=1}^{\infty}  \sum_{\iota_{1},\ldots,\iota_{s}\geq 1
} \sum_{\upsilon_{1},\ldots,\upsilon_{\tau}\geq 1} \int_{0}^{t} \int_{0}^{t_{1}-} \cdots \int_{0}^{t_{s+\tau-1-}}f_{(\upsilon_{1},\ldots,\upsilon_{\tau}, \iota_{1},\ldots,\iota_{s}, i)}(t_{1},t_{2},\ldots,t_{s+\tau })\\
 \nonumber & &\mathrm{d}G^{(\upsilon_{\tau})}_{i}(t_{s+\tau}) \ldots \mathrm{d}G^{(\upsilon_{1})}_{i}(t_{s+1})\mathrm{d}H^{(\iota_{s})}(t_{s}) \ldots \mathrm{d}H^{(\iota_{1})}(t_{1}),
\end{eqnarray*}
where $f_{(\upsilon_{1}, \ldots ,\upsilon_{\tau}, \iota_{1}, \ldots ,\iota_{s}, i)}$ is random field on $L^{2}(\Omega,\mathcal{F})$.\\
  \halmos\\

\subsection{Proof of Theorem \ref{8}}
Now, in Theorem \ref{987} we can again switch by a linear transformation from $H^{(k)}$ to $\overline{X}^{(k)}$ and from $G^{(n)}$ to $\overline{\Psi}_{i}^{(n)}$ and $\overline{\Phi}_{i}$. Thus we get the following representation:
\begin{eqnarray*}
  \lefteqn{F(t)- \mathbb{E}F(t) =\sum_{i,j=1}^{N} \sum_{s=1}^{\infty}\sum_{\tau=1}^{\infty}\sum_{\zeta=1}^{\infty}
 \sum_{\substack{(v_{1},\ldots,v_{\tau})\\  \in\{1,\ldots,b\}^{\tau}}}\sum_{\substack{(\iota_{1},\ldots,\iota_{s})\\  \in\{1,\ldots,g\}^{s}}}
  \int_{0}^{t} \int_{0}^{t_{1}-} \cdots \int_{0}^{t_{s+\tau+\zeta-1-}}}\\
 & & f_{(v_1,\ldots, v_\tau, \iota_{1},\ldots,\iota_{s},i,j)}(t, t_{1},t_{2},\ldots,t_{s+\tau+\zeta})\mathrm{d}\overline{\Phi}_{j}(t_{s+\tau+\zeta}) \ldots \mathrm{d}\overline{\Phi}_{j}(t_{s+\tau+1})\\
\nonumber & &  \mathrm{d}\overline{\Psi}^{(v_\tau)}_{i}(t_{s+\tau}) \ldots \mathrm{d}\overline{\Psi}^{(v_2)}_{i}(t_{s+2}) \mathrm{d}\overline{\Psi}^{(v_1)}_{i}(t_{s+1})\mathrm{d}\overline{X}^{(\iota_{s})}(t_{s}) \ldots \mathrm{d}\overline{X}^{(\iota_{2})}(t_{2})\mathrm{d}\overline{X}^{(\iota_{1})}(t_{1}).
\end{eqnarray*}
Using the same arguments as in the proof of Nualart and Schoutens \cite[Thm. 2]{NS} we observe that  $F- \mathbb{E}F$ can be represented
as a sum of single integrals with respect to all processes appearing in the multiple integrals.
Indeed, the procedure can be described as follows. First one takes $s=1$ and produces a sum of integrals with respect to 
$\overline{X}^{(\iota_{1})}$ for $\iota_{1}\in\{1,\ldots,g\}$. Then one takes $s=2$ and adds integrals with respect to $\overline{X}^{(\iota_{2})}$ for $\iota_{2}\in\{1,\ldots,g\}$. This procedure continues until all processes appear in the integrals.
Thus we get
\begin{eqnarray}\label{4789}
\nonumber F(t)&=&\mathbb{E}F(t)+  \int_{0}^{t}  h_{X}^{(1)}(s)\mathrm{d}\overline{X}(s)+\sum_{k=2}^{\infty}\int_{0}^{t} h^{(k)}_{X}(s)\mathrm{d}\overline{X}^{(k)}(s) +\sum_{j=1}^{N}\int_{0}^{t} h_{\Phi}^{(j)}(s)\mathrm{d}\overline{\Phi}_{j}(s)\\
 & & +\sum_{l=1}^{\infty}\sum_{i=1}^{N}\int_{0}^{t}  \tilde{h}^{(l,i)}_{\Psi}(s)\mathrm{d}\overline{\Psi}_{i}^{(l)}(s),    
\end{eqnarray}
where $h_{X}^{(1)}$, $h_{X}^{(k)}$, $h_{\Phi}^{(j)}$ and $\tilde{h}^{(l,i)}_{\Psi}$ (for $i,j=1,\ldots,N $, $k \geq 2$ and $l \geq 1$) are predictable processes. Now we want to obtain the above representation with respect to It\^{o}-Markov additive processes $X$. From the definition in (\ref{defX}) and (\ref{ovvX}), we have
\begin{equation*}
\int_{0}^{t}  h_{X}^{(1)}(s)\mathrm{d}\overline{X}(s)=\int_{0}^{t}  h_{X}^{(1)}(s)\mathrm{d}X(s)-\sum_{i=1}^{N}\int_{0}^{t}  h_{X}^{(1)}(s)\mathrm{d}\overline{\Psi}_{i}(s).
\end{equation*}
For 

\begin{equation}\label{53}
h^{(l,i)}_{\Psi}(s)= \left\{ \begin{array}{ll}
\tilde{h}^{(l,i)}_{\Psi}(s)- h_{X}^{(1)}(s)& \textrm{for $l=1$ and $i=1,\ldots,N$,}\\
\tilde{h}^{(l,i)}_{\Psi}(s)& \textrm{for $l \geq 2$ and $i=1,\ldots,N$,}
\end{array} \right.
\end{equation}
and for any square-integrable $\left\{ \mathcal{F}_t \right\}$-martingale $M$ the representation given in (\ref{4789}) can be rewritten as follows:
\begin{eqnarray}\label{b8b}
\nonumber  M(t) &=& M(0)+ \int_{0}^{t} h_{X}^{(1)}(s)\mathrm{d}X(s)+\sum_{k=2}^{\infty}\int_{0}^{t}h^{(k)}_{X}(s)\mathrm{d}\overline{X}^{(k)}(s) +\sum_{j=1}^{N}\int_{0}^{t}h_{\Phi}^{(j)}(s)\mathrm{d}\overline{\Phi}_{j}(s)\\
& & +\sum_{l=1}^{\infty}\sum_{i=1}^{N}\int_{0}^{t} h^{(l,i)}_{\Psi}(s)\mathrm{d}\overline{\Psi}_{i}^{(l)}(s),
\end{eqnarray}
where $h_{X}^{(1)}$, $h_{X}^{(k)}$, $h_{\Phi}^{(j)}$ and $h^{(l,i)}_{\Psi}$ (for $i,j=1,\ldots,N $, $k \geq 2$ and $l \geq 1$) are predictable processes.\\
 \halmos\\


\begin{thebibliography}{100}
\bibitem{App} Applebaum, D. (2004). {\it  L\'evy Processes and Stochastic Calculus}. Vol. {\bf 93} of Cambridge Studies in Advanced Mathematics. Cambridge University Press, Cambridge.
\bibitem{A2} Asmussen, S. (2003). {\it Applied Probability and Queues.}  2nd ed. Springer.
\bibitem{A1} Asmussen, S. and Albrecher H. (2010). {\it Ruin Probabilities.}  2nd ed. World Scientific, Singapore.
\bibitem{asm_avram_pist}
Asmussen, S., Avram, F. and Pistorius, M. (2004).
Russian and {A}merican put options under exponential phase-type
  {L}\'evy models. {\it Stochastic Processes and their Applications} {\bf 109}, 79--111.
\bibitem{AK} Asmussen, S. and Kella, O. (2000). Multi-dimensional martingale for Markov
additive processes and its applications. {\it Advances in  Applied Probability} {\bf 32(2)}, 376--380.
\bibitem{BNSh} Barndorff-Nielsen, O.E. and Shephard, N. (2003). Realised power variation and stochastic volatility models. {\it Bernoulli} {\bf  9}, 243--265.
\bibitem{BNSh2} Barndorff-Nielsen, O.E. and Shephard, N. (2004). {\it Financial Volatility: Stochastic Volatility and L\'evy Based Models.} Cambridge: Cambridge University Press.
\bibitem{CM} Chou, C.S. and Meyer, P.A. (1975). Sur la repr\'{e}sentation des martingales comme int\'{e}grales stochastiques dans les processus ponctuels.
{\it S\'{e}minaire de Probabilit\'{e}s} {\bf IX}, Lecture Notes in Math. 465, Springer.
\bibitem{Ca} \c Cinlar, E. (1972). Markov additive processes: I. {\it Z. Wahrscheinlichkeitstheorie und verwandte Gebiete}
{\bf 24}, 85--93.
\bibitem{Cb} \c Cinlar, E. (1972). Markov additive processes: II. {\it Z. Wahrscheinlichkeitstheorie und verwandte Gebiete}
{\bf 24}, 95-121.
\bibitem{CNS} Corcuera, J.M., Nualart, D. and Schoutens, W.  (2003).
Completion of a L\'evy market by power-jumpassets.
{\it Finance and Stochastics} {\bf 9(1)}, 109--127.
\bibitem{Da} Davis, M.H.A. (2005). The representation of martingales of jump processes.
{\it SIAM Journal on Control and Optimization} {\bf 14}, 623--638.
\bibitem{D} Davis, M.H.A. (2005). Martingale representation and all that. Systems and Control: Foundations and Applications.
In {\it Advances in Control, Communication Networks, and Transportation Systems: In Honor of Pravin Varaiya}. Birkhauser.
\bibitem{Dell1992}
Dellacherie, C., Maisonneuve, B. and Meyer, P. A. (1992). {\it Probabilit\'es et Potentiel.} Herman, Paris.
\bibitem{DOP} Di Nunno, G., Oksendal, B. and Proske, F. (2009). {\it Malliavin Calculus for L\'evy
Processes with Applications to Finance.} Springer.
\bibitem{EPQ} El Karoui, N. Peng, S. and Quenez, M.C. (1997).
Backward stochastic differential equations in finance. {\it
Mathematical Finance} {\bf 7(1)}, 1-71.
\bibitem{E1} Elliott, R.J. (1976). Stochastic integrals for martingales of a jump process with partially accessible jump times.
{\it Z. Wahrscheinlichkeitstheorie und verwandte Gebiete} {\bf 36}, 213--226.
\bibitem{E2} Elliott, R.J. (1977).
Innovation projections of a jump process and local martingales. {\it Mathematical Proceedings of the Cambridge Philosophical Society} {\bf 81}, 77--90.
\bibitem{EAM} Elliott, R.J., Aggoun, L. and Moore, J.B. (2004). {\it Hidden Markov Models: Estimation and Control.}
Springer-Verlag, Berlin.
\bibitem{emery1989}
Emery, M. (1989). On the Azema Martingales. {\it Lecture Notes in Mathematics} {\bf 1372}, 66--87, Springer, Berlin.
\bibitem{FS} F\"ollmer, H. and Schied, A. (2002). Stochastic Finance: An introduction in discrete time. {\it Studies in Mathematics} {\bf 27}. de Gruyter, Berlin-New York.
\bibitem{Palmowskietal} Gautam, N. Kulkarni, V., Palmowski, Z. and Rolski, T. (1999).
Bounds for fluid models driven by semi-Markov inputs.
{\it Probability in Engineering and Information Sciences} {\bf 13(4)}, 429--475.
\bibitem{JS} Jacod, J. and Shiryaev, A. N. (2003).
{\it Limit Theorems for Stochastic Processes.} 2nd ed. Springer-Verlag.
\bibitem{KS} Karatzas, I. and Shreve, S. (1998). {\it Brownian Motion and Stochastic Calculus}. Springer, Berlin Heidelberg New York.
\bibitem{KW} Kunita, H. and Watanabe, S. (1967). On square-integrable martingales. {\it Nagoya Mathematical Journal} {\bf 30}, 209--245.
\bibitem{Liao} Liao, M. (2004). {\it L\'evy Processes in Lie Groups.} Vol. {\bf 162} of Cambridge Tracts in Mathematics. Cambridge University Press, Cambridge.
\bibitem{LS} Liptser, R. S. and Shiryaev, A. N. (2001). {\it Statistics of Random Processes: 2. Applications.}
Springer-Verlag.
\bibitem{NS} Nualart, D. and Schoutens, W. (2000). Chaotic and predictable representations for L\'{e}vy processes.
{\it Stochastic Processes and Their Applications} {\bf 90(1)}, 109--122.
\bibitem{OS} Oksendal, B. and Sulem, A. (2004). {\it Applied Stochastic Control of Jump Diffusions.} Springer.
\bibitem{PP} Pacheco, A. and Prabhu, N.U. (1995). Markov-additive processes of arrivals.
In: Dshalalow, J.D. (ed.): {\it  Advances in Queueing}, CRC Press.
\bibitem{PTP} Pacheco, A., Tang, L.C. and Prabhu, N.U. (2009). {\it Markov-modulated Processes and Semigenerative
Phenomena.} World Scientific: Hackensack, New Jersey.
\bibitem{meiv} Ivanovs, J. and Palmowski, Z. (2012).
Occupation densities in solving exit problems for Markov additive processes and their reflections.
{\it Stochastic Processes and their Applications} {\bf 122(9)}, 3342--3360.
\bibitem{PR} Palmowski, Z. and Rolski, T. (2002). A technique for the exponential change of measure for Markov processes. {\it Bernoulli} {\bf 8(6)}, 767--785.
\bibitem{prabhu} Prabhu, N. U. (1998). {\it Stochastic Storage Processes: Queues, Insurance Risk, Dams, and
  Data Communication}. 2nd ed. Springer-Verlag.
\bibitem{P} Protter, P.E. (2005). {\it Stochastic Integration and Differential Equations: A New Approach.} 2nd ed.
Springer-Verlag, Berlin.
\bibitem{Rogers} Rogers, L. C. G. and Williams D. (2000). {\it Diffusions, Markov Processes and Martingales.}
 2nd ed. Vol. {\bf 2}. Cambridge Mathematical Library.
\bibitem{S} Schoutens, W. (1999). {\it Stochastic Processes and Orthogonal Polynomials.} Lecture Notes in Statistics,
Vol. {\bf 146}. Springer, New York.
\bibitem{ZESG} Zhang, X., Elliott, R. J., Siu, T. K. and Guo, J. Y. (2012).  Markovian regime-switching market completion using additional markov jump assets. {\it IMA Journal of Management Mathematics} {\bf 23(3)}, 283--305.
\end{thebibliography}
\end{document}